\documentclass[11pt]{article}
\usepackage{a4wide}
\usepackage{amsmath,amssymb,amsthm}
\newcommand{\R}{\mathbb{R}}
\newcommand{\N}{\mathbb{N}}
\providecommand{\PP}{\mathcal{P}}
\providecommand{\QQ}{\mathcal{Q}}
\providecommand{\E}{\mathrm{e}}
\providecommand{\D}{\mathrm{d}}
\providecommand{\esp}{\mathop{\mathrm{E}{}}\nolimits}
\providecommand{\prob}{\mathop{\mathrm{P}{}}\nolimits}
\providecommand{\abs}[1]{\lvert #1 \rvert}
\providecommand{\norm}[1]{\lVert #1 \rVert}
\providecommand{\sca}[1]{\langle #1 \rangle}

\providecommand{\tq}{\,\vert\,}

\providecommand{\card}{\mathop\mathrm{card}\nolimits}

\newtheorem{theorem}{Theorem}
\newtheorem{corollary}[theorem]{Corollary}
\newtheorem{lemma}[theorem]{Lemma}
\theoremstyle{definition}
\newtheorem{definition}[theorem]{Definition}
\title{Moments of the Gaussian Chaos}
\author{Joseph Lehec\footnote{
CEREMADE (UMR CNRS 7534)
Universit\'e Paris-Dauphine. \texttt{lehec@ceremade.dauphine.fr}
}
}
\date{}
\begin{document}

\maketitle

\begin{abstract}
This paper deals with Lata{\l}a's estimation of the moments of  Gaussian chaoses. It is shown that his argument can be simplified significantly using Talagrand's generic chaining. 
\smallskip

\noindent
Published in \emph{S\'eminaire de Probabilit\'es XLIII}, Lecture Notes in Math. 2006, Springer, 2011.
\end{abstract}

\section{Introduction}

In the article \cite{latala}, Lata{\l}a obtains an upper bound on the moments of the Gaussian chaos
\[  Y = \sum a_{n_1,\dotsc,n_d} g_{n_1} \dotsb g_{n_d} , \]
where $g_1,g_2,\dotsc$ is a sequence of independant standard Gaussian random variables and the $a_{n_1,\dotsc,n_d}$ are real numbers. His bound his sharp up to constants depending only on the order $d$ of the chaos. The purpose of the present paper is to give another proof of Lata{\l}a's result. \\
Observe that the case $d=1$ is easy, since 
\[ \bigl( \abs{\sum a_i g_i}^p \bigr)^{1/p}  =  (\sum a_i^2)^{1/2} \bigl( \esp \abs{ g_1}^p \bigr)^{1/p}
 \sim \sqrt{p} (\sum a_i^2)^{1/2} . \]
When $d=2$, Lata{\l}a recovers a result by Hanson and Wright~\cite{hanson-wright} which involves the operator and the Hilbert-Schmidt norms of the matrix $a=(a_{ij})$ 
\[ \bigl( \esp \abs{\sum a_{ij} g_i g_j}^p \bigr)^{1/p} \sim
\sqrt{p} \norm{a}_{\mathrm{HS}} + p \norm{a}_{\mathrm{op}} . \]
It is known (see \cite{pena-montgo}) that the moments of the \emph{decoupled} chaos
\[ \tilde{Y} = \sum a_{n_1,\dotsc,n_d}  g_{n_1,1} \dotsb g_{n_d,d} \]
where $( g_{i,j} )$ is a family of standard independant Gaussian variables, are comparable to those of $Y$ wih constants depending only on $d$. Using this fact and reasonning by induction on the order $d$ of the chaos, Lata{\l}a shows that the problem boils down to the estimation of the supremum of a complicated Gaussian process. Given a set $T$ and a Gaussian process $(X_t)_{t\in T}$, estimating $\esp \sup_T X_t$ amounts to studying the metric space $(T,\D)$ where $\D$ is given by the formula
\[ \D (s,t) = \bigl( \esp (X_s-X_t)^2 \bigr)^{1/2} . \] 
Dudley's estimate for instance, asserts that if the process is centered (meaning that $\esp X_t = 0$ for all $t\in T$) then there exists a universal constant $C$ such that
\[ \esp \sup T_{t} \leq C \int_0^\infty 
\sqrt{ \log N( T,\D,\epsilon ) } \ \D \epsilon ,\]
where the entropy number $N( T,\D,\epsilon )$ is the smallest number of balls (for the distance $\D$) of radius $\epsilon$ needed to cover $T$. Let us refer to Fernique~\cite{fernique} for a proof of this inequality and several applications. However, Dudley's inequality is not sharp: there exist Gaussian processes for which the integral is much larger than the expectation of the sup. Unfortunately, the phenomenon occurs here. Lata{\l}a is able to give precise bounds for the entropy numbers, but Dudley's integral does not give the correct order of magnitude. Something finer is needed.

The precise estimate of the supremum of a Gaussian process in terms of metric entropy was found by Talagrand. This was the famous \emph{Majorizing Measure Theorem}~\cite{tala87}, which is now called \emph{Generic chaining}, see the book~\cite{livre-talagrand}. Lata{\l}a did not manage to use Talagrand's theory, and his proof contains a lot of tricky entropy estimates to beat the Dudley bound. We find this part of his paper very hard to read, and our purpose is to short-circuit it using Talagrand's generic chaining.

Lastly, let us mention that we disagree with P.~Major who released an article on arXiv\footnote{{\tt http://arxiv.org/abs/0803.1453}} in which he claims that Lata{\l}a's proof is incorrect. The present paper is all about understanding Lata{\l}a's work, not correcting it.
%
\section{Notations, statement of Lata{\l}a's result}
\subsection{Tensor products, mixed injective and $L_2$ norms}
To avoid heavy multi-indices notations, it is convenient to use tensor products. 
If $X$ and $Y$ are finite dimensional normed spaces, the notation $X\otimes^\epsilon Y$ stands
for the injective tensor product of $X$ and $Y$, so that $X \otimes^\epsilon Y$ is isometric to
$\mathcal{L}(X^*,Y)$ equipped with the operator norm. If $X$ and $Y$ are Euclidean spaces, 
we denote by $X\otimes^2 Y$ their Euclidean tensor product. Moreover, in this case we identify $X$ and $X^*$, 
so that $X\otimes^2 Y$ is isometric to $\mathcal{L}(X,Y)$ equipped with the Hilbert-Schmidt norm.\\
Throughout the article $[d]$ denotes the set $\{1,\dots,d\}$.
Let $E_1,\dotsc,E_d$ be Euclidean spaces.
Given a non-empty subset $I = \{ i_1, \dots , i_p \}$ of $[d]$, we let 
\[ E_I = E_{i_1} \otimes^2 \dots \otimes^2 E_{i_p} . \] 
Also, by convention $E_\emptyset = \R$. The notation $\norm{ \cdot }_I$ stands for the norm of  $E_I$ and
\[ B_I = \{ x \in E_I;\ \norm{x}_I \leq 1 \} \] 
for its unit ball.
Let $A\in E_{[d]}$ and $\PP=\{ I_1,\dotsc,I_k \}$ be a partition of $[d]$, we let
$\norm{A}_\PP$ be the norm of $A$ as an element of the space
\[ E_{I_1} \otimes^\epsilon  \dots \otimes^\epsilon E_{I_k} . \] 
When $d=2$ for instance, the tensor $A$ can be seen as a linear map from $E_1$ to $E_2$, then
$\norm{A}_{ \{1\} \{2\} }$ and $\norm{A}_{ \{1,2\} }$ are the operator and Hilbert-Schmidt norms
of $A$, respectively. Let us give another example: assume that $d=3$ and that 
$E_1=E_2=E_3=L_2(\mu)$ for some measure $\mu$. Then
for any $f\in E_1\otimes E_2 \otimes E_3$ which we identify $L_2(\mu^{\otimes 3})$, we have
\[ \norm{f}_{\{1\}\{2,3\}} = \sup \bigl( \int f(x,y,z)  u (x) v(y,z) \ \D\mu(x)\D\mu(y)\D\mu(z) \bigr),\]
where the sup is taken over all $u,v$ having $L_2$ norms at most $1$.
Going back to the general setting, let us define for a non-empty subset $I$ of $[d]$ and an element $x\in E_I$ 
the contraction $\sca{A,x}$ to be the image of $x$ by $A$, when $A$ is seen as an element of
$\mathcal{L}(E_I,E_{[d]\backslash I})$. Then for every partition $\PP=\{ I_1,\dotsc,I_k\}$ we have
\[ \norm{A}_\PP  =
\sup \bigl\{  \sca{ A , x_1\otimes \dotsb \otimes x_k}; \ x_j\in B_{I_j} \bigr\} . \]
If $\QQ=\{J_1,\dotsc, J_l\}$ is a finer partition than $\PP$ 
(this means that any element of $\QQ$ is contained in an element of $\PP$) then
\[ \{ x_1\otimes \dots \otimes x_l , \  x_j\in B_{J_j} \} \subset 
\{ y_1\otimes \dots \otimes y_k , \  y_j\in B_{I_j}\} , \]
hence $\norm{A}_\QQ \leq \norm{A}_\PP$.
In particular,
\[ \norm{A}_{\{1\}\dotsb\{d\}} \leq \norm{A}_\PP \leq \norm{A}_{[d]} . \]
\subsection{Moments of the Gaussian chaos}
If $\PP$ is a partition of $[d]$, its cardinality $\card \PP$ is the number of
subsets of $[d]$ in $\PP$. Let $E_1,\dotsc,E_d$ be Euclidean
spaces and $A\in E_{[d]}$. Let $X_1,\dotsc,X_d$ be independant random vectors
such that for all $i$, the vector $X_i$ is a standard Gaussian vector of $E_i$.
The (real) random variable 
\[ Z =  \sca{A , X_1 \otimes \dotsb \otimes X_d }  \]  
is called \emph{decoupled} Gaussian chaos of order $d$. Here is the main result of Lata{\l}a.
\begin{theorem}
\label{1main}
There exists a constant $\alpha_d$ depending only on $d$ such that for all $p\geq 1$
\[ \bigl( \esp \abs{ Z } ^p \bigr)^{1/p}
   \leq  \alpha_d \sum_\PP p^{ \frac{\card \PP}{2}}  
  \norm{ A }_\PP  , \]
the sum running over all partitions $\PP$ of $[d]$.
\end{theorem}
The following theorem and corollary are intermediate results from which the previous theorem shall follow; however we believe they are of independent interest.
\begin{theorem}
\label{1bis}
Let $F_1,\dotsc,F_{k+1}$ be Euclidean spaces,
let $A\in F_{[k+1]}$ and $X$ be a standard Gaussian vector on $F_{k+1}$, recall that $\sca{A,X}\in F_1\otimes \dotsb \otimes F_k$.
Then for all $\tau \in (0,1)$:
\[  \esp \norm{ \sca{  A , X  } }_{\{1\}\dotsb\{k\}}  \leq  \beta_k \sum_\PP
                \tau^{k-\card \PP}  \norm{A}_\PP  , \]
where the sum runs over all partitions $\PP$ of $[k+1]$ and the constant $\beta_k$ depends only on $k$.
\end{theorem}
\begin{corollary}
\label{cor}
Under the same hypothesis, we have for all $p\geq 1$
\[ \Bigl( \esp \norm{ \sca{ A , X } }^p_{\{1\}\dotsb\{k\}} \Bigr)^{1/p} \leq  \delta_k
\sum_{\PP} p^{\frac{\card \PP-k}{2}} \norm{A}_\PP . \]
\end{corollary}
\begin{proof}
Let $f\colon x\in F_{k+1} \mapsto \norm{ \sca{A,x} }_{\{1\}\dotsb\{k\}}$. Let us
use the concentration property of the Gaussian measure, which asserts that Lipschitz
functions are close to their means with high probability. More 
precisely, letting $m=\esp f(X)$, we have for all $p\geq 1$
\[ 
\bigl( \esp \abs{f(X)-m}^p \bigr)^{1/p} \leq
    \delta' \sqrt{p} \norm{f}_{\mathrm{lip}}, 
\]
where $\norm{f}_{\mathrm{lip}}$ is the Lipschitz constant of $f$ and $\delta'$
is a universal constant. We refer to \cite{ledoux} for more details on this inequality. 
Noting that
\[ \norm{f}_{\mathrm{lip}} = \sup_{ x\in B_{k+1} } \  \norm{\sca{A,x}}_{\{1\}\dotsb\{k\}}
  = \norm{A}_{\{1\}\dotsb\{k+1\}} . \]
and using the triangle inequality, we get
\[  \bigl( \esp \abs{f(X)}^p \bigr)^{1/p} \leq
     \esp f(X) +  \delta' \sqrt{p} \norm{A}_{\{1\}\dotsb\{k+1\} } . \]
The result then follows from the upper bound on $\esp f(X)$ given by Theorem~\ref{1bis} with $\tau=p^{-1/2}$.
\end{proof}
Let us prove Theorem~\ref{1main}. We proceed by induction on $d$.
When $d=1$, the random variable $\sca{A,X_1}$ is, in law, equal to the Gaussian variable of 
variance $\norm{A}_{\{1\}}^2$. The $p$-th moment of the standard Gaussian variable being of order $\sqrt{p}$, we get
\[
\bigl( \esp \abs{ \sca{ A, X_1 } }^p \bigr)^{1/p} \leq \alpha \sqrt{p} \norm{A}_{\{1\}} \]
for some universal $\alpha$, hence the theorem for $d=1$.
\\
Assume that the result holds for chaoses of order $d-1$. From now on, if $I=\{i_1,\dotsc,i_r\}$ is a subset of
$[d]$ we denote the tensor $X_{i_1} \otimes \dotsb \otimes X_{i_r}$ by $X_I$. Notice that
\[ \sca{A , X_{[d]} } = \bigl\langle  \sca{A, X_d} , X_{[d-1]} \bigr\rangle \]
and apply the induction assumption to the matrix $B= \sca{A,X_d}$. This yields
\[
\esp \bigl( \abs{ \sca{B , X_{[d-1]}} }^p \big\vert X_d \bigr) \leq \alpha_{d-1}^p
    \Bigl( \sum_\PP p^{ \frac{\card \PP}{2} }
       \norm{B}_\PP \Bigr)^p  ,
\]
where the sum runs over all partitions $\PP$ of $[d-1]$. Taking expectation and the $p$-th root, we obtain
\begin{equation}
\label{1recurrence}
\begin{split}
\bigl( \esp \abs{  \sca{A, X_{[d]}}  }^p \bigr)^{1/p} & \leq
\alpha_{d-1} \left( 
   \esp \Bigl( \sum_\PP p^{ \frac{\card \PP}{2} } 
         \norm{ \sca{A, X_d} }_\PP \Bigr)^p \right)^{1/p}
 \\
& \leq  \alpha_{d-1} \sum_\PP p^{ \frac{\card \PP}{2} }
  \Bigl( \esp \norm{ \sca{A, X_d} }^p_\PP \Bigr)^{1/p} ,
\end{split}
\end{equation}
by the triangle inequality.
Let $\PP=\{I_1,\dotsc,I_k\}$ be a partition of $[d-1]$.
Let $F_i = E_{I_i}$ for $i\in[k]$ and $F_{k+1}=E_{d}$.
The tensor $A$ can be seen as an element of $F_{[k+1]}$, let us rename it
$A'$ when we do so. Corollary~\ref{cor} gives
\[ 
\bigl( \esp \norm{ \sca{ A' , X_d} }^p_{\{1\}\dotsb\{k\}} \bigr)^{1/p}  \leq
    \delta_k p^{-\frac{k}{2}} \sum_\QQ p^{ \frac{\card \QQ}{2} } \norm{ A' }_{\QQ} ,
\]
where the sum is taken over all partitions $\QQ$ of $[k]$.
Going back to the the space $E_{[d]}$, this inequality translates as
\begin{equation}
\label{1etapecor}
\bigl( \esp \norm{ \sca{A , X_d} }^p_\PP \bigr)^{1/p}  \leq
    \delta_k p^{-\frac{k}{2}} \sum_\QQ p^{ \frac{\card \QQ}{2} } \norm{A}_{\QQ} ,
\end{equation}
and this time the sum runs over partitions $\QQ$ of $[d]$ such that the partition 
\[ \bigl\{ I_1,\dotsc,I_k,\{d\} \bigr\} \]
is finer than $\QQ$. However, the inequality still holds if we 
take the sum over all partitions of $[d]$ instead. 
We plug \eqref{1etapecor} into \eqref{1recurrence}, 
the numbers $p^{\frac{\card\PP}{2}}$ cancel out and we get the desired inequality
with constant
\[  \alpha_d = \alpha_{d-1} \sum_\PP \delta_{\card \PP} , \]
where the sum is taken over all partitions $\PP$ of $[d-1]$.

So it is enough to prove Theorem~\ref{1bis}, this is the purpose of the rest of the article.
\section{The generic chaining}
Let $F_1,\dotsc,F_{k+1}$ be Euclidean spaces, let $A\in F_{[k+1]}$ and $X$ be a standard Gaussian vector of $F_{k+1}$. For $i\in [k]$ let $B_i$ be the unit ball of $F_i$, let $T=B_1 \times \dotsb \times B_k$. Recall that for $x=(x_1,\dotsc,x_k) \in T$, the notation $x_{[k]}$ stands for the tensor $x_1\otimes \dotsb \otimes x_k$. Note that
\begin{equation}
\label{1process1}
\esp \norm{ \sca{ A , X} }_{\{1\}\dotsb\{k\}} = \esp \sup_{x\in T} \sca{ A , x_{[k]}\otimes X } = 
\esp \sup_{x\in T} \bigl\langle \sca{ A,x_{[k]} }  , X\bigr\rangle  . 
\end{equation}
Notice also that $(P_x)_{x\in T} = \bigl( \bigl\langle \sca{ A,x_{[k]} }  , X\bigr\rangle \bigr)_{x\in T}$ is a Gaussian process. 
To estimate $\esp \sup_T P_x$, we shall study the metric space $(T,\D)$, where 
\[ \D(x,y) = \bigl( \esp (P_x-P_y)^2 \bigr)^{1/2} . \]
This distance can be computed explicitly. Indeed
\begin{equation}
\label{1process2}
\D (x,y)^2 =  
\esp \bigl\langle \sca{ A,x_{[k]} - y_{[k]} }  , X \bigr\rangle^2 =
\norm{ \sca{ A , x_{[k]}-y_{[k]} } }_{ \{ k+1 \} }^2 .  
\end{equation}
The \emph{generic chaining}, introduced by Talagrand, will be our main tool. We sketch briefly the main ideas of the theory and refer to Talagrand's book~\cite{livre-talagrand} for details.\\
Let $(T,\D)$ be a metric space. If $S$ is a subset of $T$ we let $\delta_\D(S)$ be the diameter of $S$
\[ \delta_\D(S) = \sup_{s,t\in S} \ \D(s,t) . \]
Given a sequence $\bigl( \mathcal{A}_n \bigr)_{n\in \N}$ of partitions of $T$ and an element $t\in T$, we let $A_n (t)$ be the unique element of $\mathcal{A}_n$ containing $t$.
\begin{definition}
Let
\[ \gamma_\D (T) = \inf \Bigl( \sup_{t\in T}  \sum_{n=0} ^\infty
             \delta_\D \bigl( A_n(t) \bigr) 2^{n/2} \Bigr) , \]
where the infimum is over all sequences of partitions $\bigl( \mathcal{A}_n \bigr)_{n\in \N}$ of $T$ satisfying the cardinality condition
\begin{equation}
\label{cardinality}
\mathcal{A}_0 = \{ T \} \quad \mathrm{and} \quad \forall n\geq 1, \  \card \mathcal{A}_n \leq 2^{2^n} .
\end{equation}
\end{definition}
Notice that $\gamma_\D(T) \geq \delta_\D(T)$. In particular, if the metric is not trivial then $\gamma_\D (T)$ is non-zero. Thus there exists a sequence of partitions $(\mathcal{A}_n)_{n\in\N}$ satisfying the cardinality condition and
\[ \sup_{t\in T} \sum_{n=0} ^\infty
             \delta_\D \bigl( A_n(t) \bigr) 2^{n/2}  \leq 2 \gamma_\D (T) . \]
We recall the all important
\begin{theorem}[Majorizing Measure]
There exists a universal constant $\kappa$ such that for any Gaussian process $(X_t)_{t\in T}$ that is centered (meaning $\esp X_t = 0$ for all $t\in T$) we have 
\[ \tfrac{1}{\kappa} \gamma_{\D}  (T)   \leq  
  \esp \sup_{t\in T} X_t 
\leq \kappa \gamma_{\D} (T) , \]
where the metric $\D$ is defined by $\D (s,t) = \bigl( \esp (X_s-X_t)^2 \bigr)^{1/2}$.
\end{theorem}
Here are two simple lemmas.
\begin{lemma}
\label{1reduction}
Let $(T,\D)$ be a metric space. Let $a,b \geq 1$, and 
$(\mathcal{A}_n)_{n\in \N}$ be a sequence of partitions of $T$ satisfying
\[ \forall n\in \N, \  \card \mathcal{A}_n \leq 2^{a + b 2^n} . \] 
Letting $\gamma =\sup_{t\in T} \sum_{n=0}^\infty 
     \delta_\D  \bigl( A_n (t) \bigr) 2^{n/2}$, we have
\[ \gamma_\D (T) \leq 
        \rho \bigl( \sqrt{ab} \ \delta_\D (T) + \sqrt{b} \ \gamma \bigr) , \]
for some universal $\rho$.
\end{lemma}
\begin{proof}
Let $p,q$ be the smallest integers satisfying $a\leq 2^p$ and $b\leq 2^q$. Let
\[
\mathcal{B}_n = 
\left\{
\begin{array}{ll}
\{ T \}    & \text{ if } n\leq p+q \\
\mathcal{A}_{n-q-1} & \text{ if } n\geq p+q+1. 
\end{array} \right.
\]
If $n\geq p+q+1$ then $p \leq n-1$ so 
\[ \card \mathcal{B}_n \leq 2^{2^p + 2^{n-1}} \leq 2^{2^n} . \]  
Thus the  sequence $(\mathcal{B}_n)_{n\in\N}$ satisfies \eqref{cardinality}. On the other hand, for all $t\in T$
\begin{align*} 
\sum_{n=0}^\infty \delta_\D \bigl( B_n (t) \bigr) 2^{n/2}
  & =  \sum_{n=0}^{p+q} \delta_d (T)  2^{n/2} + \sum_{n=p}^\infty 
     \delta_\D \bigl( A_n (t) \bigr)  2^{\frac{n+q+1}{2}}  \\
  & \leq \tfrac{ 2^{\frac{p+q+1}{2}} -1 }{\sqrt{2}-1} \delta_\D (T) 
         +  2^{\frac{q+1}{2}}  \gamma .
\end{align*}
Moreover $2^p\leq 2a$ and $2^q\leq 2b$, hence the result.
\end{proof}
\begin{lemma}
\label{1somme}
Let $\D_1,\dots,\D_N$ be distances defined on $T$ and let $\D=\sum \D_i$. Then
\[ \gamma_{\D} (T) 
    \leq \rho' \sqrt{N}  \sum_{i=1}^N \gamma_{\D_i} (T) ,\]
where $\rho'$ is a universal constant.
\end{lemma}
\begin{proof}
For all $i \in[N]$, there exists a sequence $(\mathcal{A}_n^i)_{n\in \N}$ of partitions of $T$ satisfying the cardinality condition~\eqref{cardinality} and
\[ \sup_{t\in T} \sum_{n=0}^\infty \delta_{\D_i} \bigl( A_n^i(t) \bigr) 2^{n/2}  \leq 2 \gamma_{\D_i} ( T ). \]
Then let
\[ \mathcal{A}_n = \{ A^1\cap\dotsb \cap A^N, \  
                          A^i \in \mathcal{A}_n^i \} . \]
This clearly defines a sequence of partitions of $T$, and for all $n$ we have
\begin{equation}   
\label{1card0}
\card \mathcal{A}_n \leq  2^{ N 2^n } .
\end{equation}
On the other hand, for all $t\in T$ and $i\in [N]$ we have
$A_n(t) \subset A_n^{i} (t)$, so 
\[ \delta_\D \bigl( A_n(t) \bigr)
 \leq   \sum_{i=1}^N \delta_{\D_i} \bigl( A_n (t) \bigr)
 \leq   \sum_{i=1}^N \delta_{\D_i} \bigl( A_n^i (t) \bigr)  . \]
Consequently
\begin{equation}
\label{1delta0}
\sup_{t\in T} \sum_{n=0}^\infty \delta_{\D}
        \bigl( A_n (t) \bigr) 2^{n/2} \leq
     2    \sum_{i=1}^N \gamma_{\D_i}( T ). 
\end{equation}
By the previous lemma, equations \eqref{1card0} and \eqref{1delta0} yield the result.
\end{proof}
\section{Proof of Theorem~\ref{1bis}}
The proof is by induction on $k$. When $k=1$ the theorem is a consequence of the following: let $A \in F_1\otimes F_2$ and $X$ be a standard Gaussian vector on $F_2$, then
\[
\esp \norm{ \sca{ A , X} }_{\{1\}}
   \leq  \bigl( \esp \norm{ \sca{ A , X } }^2_{\{1\}} \bigr)^{1/2}  
     = \norm{A}_{\{1,2\}} .
\]
Assume that $k\geq 2$ and that the theorem holds for any $k'<k$. Let $A\in F_{[k+1]}$. Recall that for $i\in [k]$ the unit ball of $F_i$ is denoted by $B_i$ and the product $B_1\times\dots\times B_k$ by $T$. Let $I$ be a \emph{non-empty} subset of $[k]$ and $\D_I$ be the pseudo-metric on $T$ defined by
\begin{equation}
\label{1distanceI}
\D_I (x,y) = \norm{ \sca{ A , x_I - y_I } }_{[k+1]\backslash I} . 
\end{equation}
By the majorizing measure theorem and the equations \eqref{1process1} and \eqref{1process2}, Theorem~\ref{1bis} is equivalent to
\begin{theorem}
\label{thmbis}
For all $\tau\in (0,1)$ 
\[ \gamma_{ \D_{[k]} } ( T ) \leq \beta'_k \sum_\PP \tau^{k-\card \PP} \norm{A}_\PP , \]
with a sum running over all partitions $\PP$ of $[k+1]$.
\end{theorem}
Our purpose is to prove Theorem~\ref{thmbis} by induction on $k$. Let $\tau$ be a fixed positive real number and let $\D^\tau$ be the following metric: 
\begin{equation}
\label{dtau}
\D^\tau = \sum_{\emptyset \subsetneq I \subsetneq [k]} \tau^{k-\card I} \D_I   .
\end{equation}
Let us sketch the argument. First we use an entropy estimate and the generic chaining to compare  $\gamma_{\D_{[k]}} (T)$ and $\gamma_{\D^\tau} (T)$, then we use the induction assumption to estimate the latter.\\
Here is the crucial entropy estimate of Lata{\l}a~\cite[Corollary 2]{latala}.
\begin{lemma}
\label{entropie}
Let $S\subset T$, let $\tau \in (0,1)$ and $\epsilon=\delta_{\D^\tau} (S)  + \tau^k \norm{A}_{[k+1]}$. Then
\[ N \bigl( S, \D_{[k]}, \epsilon \bigr) 
      \leq 2^{ c_k \tau^{-2} } , \]
for some constant $c_k$ depending only on $k$.
\end{lemma}
Let us postpone the proof to the last section. \\
Let $(\mathcal{B}_n)_{n\in \N}$ be a sequence of partitions of $T$ satisfying the cardinality condition \eqref{cardinality} and
\begin{equation}
\label{1gammaB}
\sup_{t \in T} \sum_{n=0}^\infty \delta_{ \D^\tau } \bigl( B_n(t) \bigr) 
2^{n/2} \leq 2 \gamma_{ \D^\tau }  (T) .
\end{equation}
Let $n\in \N$ and $B \in \mathcal{B}_n$, set $\tau_n= \min (\tau, 2^{-n/2})$ and $\epsilon_n = \delta_{ \D^{\tau_n} } (B) + \tau_n^k\norm{A}_{[k+1]}$. Observe that $\tau_n^{-2} \leq \tau^{-2} + 2^n$ and apply Lemma~\ref{entropie} to $B$ and $\tau_n$:
\[ N( B, \D_{[k]} , \epsilon_n ) \leq 2^{ c_k \tau_n^{-2} } \leq 2^{ c_k \tau^{-2} + c_k 2^n} . \]
Therefore we can find a partition $\mathcal{A}_B$ of $B$ whose cardinality is controlled by the number above and such that any $R\in\mathcal{A}_B$ satisfies
\[ \delta_{\D_{[k]} }  (R) \leq 2 \epsilon_n \leq  2 \delta_{ \D^\tau }(B) + 2 \tau_n^k \norm{A}_{[k+1]} . \]
Indeed $\tau_n \leq \tau$ implies that $\D^{\tau_n} \leq \D^\tau$.
Then we let $\mathcal{A}_n= \cup \{ \mathcal{A}_B ; \  B \in \mathcal{B}_n \}$.
This clearly defines a sequence of partitions of $T$ which satisfies
\begin{align} 
\label{card}
\card \mathcal{A}_n & \leq 2^{ c_k \tau^{-2} + c_k 2^n } \card \mathcal{B}_n \leq 2^{ c_k \tau^{-2} + (c_k +1) 2^n },\\
\label{diam}
\delta_{\D_{[k]} } \bigl( A_n(t) \bigr) & \leq
  2 \delta_{\D^\tau} \bigl( B_n(t) \bigr) + 2 \tau_n^k \norm{A}_{[k+1]} ,
\end{align}
for all $n\in \N$ and $t\in T$. 
Recall that $\tau_n=\min(\tau,2^{-n/2})$, an easy computation shows that
\[ 
\sum_{n=0}^\infty \tau_n^k 2^{n/2} \leq C \tau^{k-1}  
\]
for some universal $C$. Therefore, for all $t\in T$, we have
\begin{align*}
\sum_{n=0}^\infty 
    \delta_{\D_{[k]} } \bigl( A_n (t)\bigr) 2^{n/2}  & \leq
    2 \sum_{n=0}^{\infty} \bigl( \delta_{\D^\tau} ( B_n (t) ) +  \tau_n^k \norm{A}_{[k+1]} \bigr) 2^{n/2} , \\
 & \leq 4 \gamma_{\D^\tau} (T) +2 C \tau^{k-1} \norm{A}_{[k+1]}. 
\end{align*}
By \eqref{card} and applying Lemma~\ref{1reduction}, we get for some constant $C_k$ depending only on $k$
\begin{equation}
\label{gamma}
\begin{split}
\gamma_{\D_{[k]}} (T) & \leq C_k \bigl( \gamma_{\D^\tau} (T) + \tau^{k-1} \norm{A}_{[k+1]} + 
\tau^{-1} \delta_{\D_{[k]}}(T) \bigr) , \\
                    & \leq 2 C_k \bigl( \gamma_{\D^\tau}(T) + \tau^{k-1} \norm{A}_{[k+1]} + 
\tau^{-1} \norm{A}_{\{1\} \dotsb \{ k+1\}} \bigr) .
\end{split}
\end{equation}
Indeed
\[
\delta_{\D_{[k]}} (T) = 2 \sup_{x\in T} \norm{ \sca{ A, x } }_{\{ k+1\}}
    = 2 \norm{A}_{\{1\} \dotsb \{ k+1\}} .
\]
We have not used the induction assumption yet.
Let $I=\{i_1,\dotsc,i_p\}$ be a subset of $[k]$, different from $\emptyset$ and $[k]$. For $j\in[p]$ let $F'_j = F_{i_j}$ and
let $F'_{p+1} = F_{ [k+1]\backslash I }$. Since $p<k$ we can apply inductively Theorem~8 to the tensor $A$ seen as an element of $F'_{[p+1]}$. For all $\tau\in (0,1)$
\begin{equation}
\label{1gammadI}
\gamma_{\D_I} ( T ) \leq \beta'_p \sum_\QQ \tau^{p - \card \QQ} \norm{A}_{\QQ} ,
\end{equation}
where the sum runs over all partitions $\QQ$ of $[k+1]$ such that the partition $\{i_1\}, \dotsc ,\{i_p\}, [k+1]\backslash I$ is finer than $\QQ$.
Again, the inequality is still true if we take the sum over all partitions of $[k+1]$ instead. 
According to Lemma~\ref{1somme} and since $\gamma$ is clearly homogeneous, we have 
\[ \gamma_{\D^\tau} ( T )  \leq \rho' 
 \sqrt{N} \sum_{\emptyset \subsetneq I \subsetneq [k]} \tau^{k-\card I} \gamma_{\D_I} (T)  \]
where $N$ is the number of subsets of $[k]$ which are different from $\emptyset$ and $[k]$, namely $2^k-2$. By \eqref{1gammadI} we get
\[
\gamma_{\D^\tau} ( T )  \leq D_k \sum_\PP \tau^{k-\card \PP} \norm{A}_\PP , 
\]
for some $D_k$ depending only on $k$. This, together with~\eqref{gamma}, concludes the proof of Theorem~\ref{thmbis}.

In the last section we prove Lemma~\ref{entropie}, this is essentially Lata{\l}a's proof.
\section{Proof of the entropy estimate}
Let $x=(x_1,\dots,x_k)\in F_1\times\dots\times F_k$, let $\abs{x_i}$ be the norm of $x_i$ in $F_i$.  Let $X_1,\dots,X_k$ be independant standard Gaussian vectors on $F_1,\dotsc,F_k$, respectively.
\begin{lemma}
\label{volume}
For all semi-norm $\norm{ \cdot }$ on $F_{[k]}$, we have
\[ \prob \Bigl( \norm{ X_{[k]} - x_{[k]} }
           \leq \esp \sum_{\emptyset \subsetneq I \subset [k]}
       4^{\card I} \norm{ X_I \otimes x_{[k]\backslash I} } \Bigl)
   \geq 2^{-k} \E^{- \frac{1}{2} \sum_{i=1}^k \abs{x_i}^2} . \]
\end{lemma}
\begin{proof}
Let us start with an elementary remark.
Let $x\in \R^n$, let $K$ be a symmetric subset of $\R^n$, and $\gamma_n$ be the standard Gaussian measure on $\R^n$. Then
\begin{equation}
\label{translate}
\gamma_n ( x+K ) \geq \gamma_n (K) \E^{- \frac{1}{2} \abs{x}^2 } .
\end{equation}
Indeed, the symmetry of $K$ the convexity of the exponential function imply that
\begin{align*}  
\int_{x+K} \E^{-\frac{1}{2} \abs{z}^2} \, dz & = 
\int_K \tfrac{1}{2}( \E^{ - \frac{1}{2} \abs{x+y}^2}
                  + \E^{ - \frac{1}{2} \abs{x-y}^2} ) \, dy  \\
& \geq \int_K \E^{ - \frac{1}{2} ( \abs{x}^2 + \abs{y}^2 ) }  \, dy 
\end{align*}
which proves \eqref{translate}. \\
Let us prove the lemma by induction on $k$. If $k=1$, applying \eqref{translate} to 
$K=\{y\in F_1, \  \norm{y} \leq 4 \esp \norm{X_1}\}$ and $x=x_1$, we get
\[ 
\prob \bigl( \norm{ X_1-x_1 } \leq 4 \esp \norm{X_1} \bigr) \geq 
      \E^{-\frac{1}{2} \abs{x_1}^2} \prob \bigl( \norm{ X_1 } \leq 4 \esp \norm{X_1} \bigr).
\]
Besides, by Markov we have
$\prob \bigl( \norm{X_1} \geq 4 \esp \norm{X_1} \bigr) \leq \frac{1}{4} \leq \frac{1}{2}$,
hence the result for $k=1$. \\
Let $k\geq 2$ and assume that the result holds for $k-1$. Let 
\begin{align*}
S & = \sum_{\emptyset\subsetneq I\subset [k-1]} 4^{\card I} 
   \norm{ X_{I} \otimes x_{[k-1]\backslash I} \otimes X_k} \\
T & = \sum_{\emptyset\subsetneq I\subset [k-1]} 4^{\card I} 
   \norm{ X_{I} \otimes x_{[k-1]\backslash I} \otimes x_k } 
\end{align*}
and let $A$, $B$ and $C$ be the events
\begin{align*}
A & = \bigl\{ \norm{ x_{[k-1]} \otimes (X_k-x_k) } \leq 
           4 \esp \norm{ x_{[k-1]} \otimes X_k } \bigr\} \\
B & = \bigl\{ \norm{ ( X_{[k-1]} - x_{[k-1]} ) \otimes X_k } \leq 
  \esp ( S \, \vert \, X_k ) \bigr\}\\
C & = \bigl\{ \esp ( S \, \vert \, X_k ) \leq 4 \esp S + \esp T \bigr\}.
\end{align*}
By the following triangle inequality
\[ \norm{ X_{[k]} - x_{[k]} } \leq \norm{ x_{[k-1]} \otimes(X_k-x_k) } + \norm{ (X_{[k-1]}-x_{[k-1]})\otimes X_k } , \]
when $A$, $B$ and $C$ occur we have
\begin{align*}
\norm{ X_{[k]} - x_{[k]} } & \leq 4 \esp \norm{ x_{[k-1]} \otimes X_k } + 4 \esp S + \esp T \\
                          & = \esp \sum_{\emptyset \subsetneq I \subset [k]} 4^{\card I}
                   \norm{ X_I \otimes x_{[k]\backslash I} } .
\end{align*}
Assume that $X_k$ is deterministic, and apply the induction assumption 
to the spaces $F_1,\dotsc,F_{k-1}$ and to the semi-norm $\norm{ y }_1 = \norm{ y\otimes X_k}$ for all $y\in F_{[k-1]}$, then
\[ \prob( B \tq X_k ) \geq 2^{-k+1}
\E^{-\frac{1}{2}\sum_{i=1}^{k-1} \abs{x_i}^2} .\] 
Since $A$ and $C$ depend only on $X_k$, this implies that
\[ \prob( A \cap B \cap C ) \geq \prob ( A \cap C) 2^{-k+1}
\E^{-\frac{1}{2}\sum_{i=1}^{k-1} \abs{x_i}^2} . \]
So it is enough to prove that $\prob ( A \cap C) \geq 2^{-1} \E^{- \frac{1}{2} \abs{x_k}^2 }$.
For all $y \in F_k$ we let
\begin{align*} 
\norm{ y}_2 & = \norm{ x_{[k-1]} \otimes y }, \\
\norm{y}_3 & =  \esp \sum_{\emptyset \subsetneq I\subset [k-1]} 4^{\card I} 
   \norm{ X_{I} \otimes x_{[k-1]\backslash I} \otimes y } .
\end{align*}
So that
\begin{align*}
A & = \bigl\{ \norm{ X_k - x_k }_2 \leq 4 \esp \norm{ X_k }_2 \bigr\} , \\
C & = \bigl\{ \norm{ X_k }_3 \leq 4 \esp \norm{ X_k }_3 + \norm{x_k}_3 \bigr\} .
\end{align*}
Let
\[ K = \{ y \in F_k, \ \norm{y}_2 \leq 4 \esp \norm{X_k}_2\} \cap 
       \{ y \in F_k, \  \norm{y}_3 \leq 4 \esp \norm{X_k}_3\} ,
\]
then, by the triangle inequality, the event $X_k \in x_k +K $ is included in $A\cap C$. Using \eqref{translate}, we get
\[ \prob (A\cap C) \geq \prob ( X_k \in x_k + K) \geq \E^{- \frac{1}{2} \abs{x_k}^2 } \prob (X_k \in K) . \] 
Therefore, it is enough to prove that $\prob(X_k \in K) \geq \frac{1}{2}$, and this is a simple application of Markov again.
\end{proof}
Let $F_{k+1}$ be another Euclidean space and let $A\in F_{[k+1]}$. Recall that for $I=\{i_1,\dotsc,i_p\} \subset [k+1]$, we let 
\[ F_I = F_{i_1} \otimes^2 \dotsb \otimes^2 F_{i_p} \]
and $\norm{ \cdot }_I$ be the corresponding (Euclidean) norm. Our purpose is to apply the previous lemma to the semi-norm defined by $\norm{y} = \norm{\sca{  A , y }}_{ \{k+1\} }$, for all $y\in F_{[k]}$. Notice that for all $x\in F_1\times \dotsb \times F_k$ and for all $\emptyset \subsetneq I\subsetneq [k]$
\begin{align*}
\esp \norm{ X_I \otimes x_{[k] \backslash I} } 
& \leq \bigl( \esp \norm{ X_I \otimes x_{[k] \backslash I} }^2 \bigr)^{1/2} \\
& = \norm{ \sca{ A, x_{[k]\backslash I} } }_{ I \cup \{k+1\} } ,
\end{align*}
which, according to the definition~\eqref{1distanceI}, is equal to $\D_{ [k]\backslash I} (0,x)$. In the same way, when $I=[k]$
\[ \esp \norm{ \sca{ A , X_{[k]} } }_{ \{k+1\} } \leq \norm{ A }_{ [k+1] } . \] 
We let the reader check that Lemma \ref{volume} then implies the following:
for all $\tau\in (0,1)$ and $x\in T$, letting $\epsilon_x = \D^\tau (x,0) + \tau^k \norm{A}_{[k+1]}$, we have 
\begin{equation}
\label{volumetric}
\prob \bigl( \D_{[k]} ( x, \tau X ) \leq \epsilon_x  / 2\bigr) \geq
      2^{ - c_k \tau^{-2}} 
\end{equation}
for some constant $c_k$ depending only on $k$.\\
Lemma \ref{entropie} follows easily from this observation. Indeed let $S\subset T$, since $S$ and its translates have the same entropy numbers, we can assume that $0\in S$. Then $\epsilon_x  \leq \epsilon := \delta_{\D^\tau} (S) + \tau^k \norm{A}_{[k+1]}$ for all $x\in S$. Let $S'$ be a subset of $S$ satisfying
\begin{itemize}
\item[(i)] For all $x,y\in S'$, $\D_{[k]} (x,y) \geq \epsilon$.
\item[(ii)] The set $S'$ is maximal (for the inclusion) with this property. 
\end{itemize}
By maximality $S'$ is an $\epsilon$-net of $S$, so $N( S, \D_{[k]}, \epsilon ) \leq \card S'$. On the other hand, by (i) the balls (for $\D_{[k]}$) of radius $\epsilon/2$ centered at different points of $S'$ do not intersect. This, together with \eqref{volumetric}, implies that 
\[ 2^{ - c_k \tau^{-2} } \card S' \leq   \sum_{x \in S'} \prob \bigl( \D_{[k]} ( x, \tau X ) \leq  \epsilon / 2\bigr) \leq 1 , \]
hence the result.

\end{document}